\documentclass[12pt]{article}
\usepackage{amssymb,oldlfont}
\def\wrtext#1{\relax\ifmmode{\leavevmode\hbox{#1}}\else{#1}\fi}
\def\abs#1{\left|#1\right|}
\def\begeq{\begin{equation}}
\def\endeq{\end{equation}}
\def\Remark{\vskip 2mm \noindent {\em Remark}}

\newcommand{\eps}{\varepsilon}
\def\part#1{\frac{\partial}{\partial #1}}

\def\norm#1{||\,#1\,||}

\newcommand{\hk}[1]{\langle #1\rangle}

\newcommand{\real}{\mbox{\bf R}}
\newcommand{\comp}{\mbox{\bf C}}
\newcommand{\z}{\mbox{\bf Z}}
\newcommand{\nat}{\mbox{\bf N}}

\newcommand{\dist}{\mbox{\rm dist\,}}
\newcommand{\Spec}{\mbox{\rm Spec\,}}
\newcommand{\mod}{\mbox{\rm mod\,}}
\renewcommand{\Re}{\mbox{\rm Re\,}}
\renewcommand{\Im}{\mbox{\rm Im\,}}
\renewcommand{\exp}{\mbox{\rm exp\,}}
\newcommand{\supp}{\mbox{\rm supp}}

\def\neigh{neighborhood}

\def\Re{{\rm Re\,}}
\def\Im{{\rm Im\,}}

\newtheorem{dref}{Definition}[section]

\newtheorem{theo}[dref]{Theorem}
\newtheorem{prop}[dref]{Proposition}

\newenvironment{proof}{\vspace{.3cm}\noindent{{\em Proof:}}}{\hfill$\Box$
\vspace{.2cm}}

\begin{document}
\begin{center}
\Large {\bf Eigenfrequencies and expansions for damped wave equations}
\end{center}
\vspace*{3mm}
\noindent
{\bf Michael Hitrik}

\vspace*{3mm}
\noindent
Department of Mathematics, University of California, Los Angeles, 
CA 90095-1555, USA
\vskip 1mm
\noindent

\vspace*{1cm}
\noindent
{\bf Abstract}: We study eigenfrequencies and propagator expansions
for damped wave equations on compact manifolds. Under the assumption 
of geometric control, the propagator is shown to
admit an expansion in terms of finitely many eigenmodes near the real axis, with an error term
exponentially decaying in time. In the presence of a nondegenerate elliptic closed
geodesic not meeting the support of the damping coefficient, we show
that there exists a sequence of eigenfrequencies converging rapidly to
the real axis. In the case of Zoll manifolds, we show that the propagator can
be expanded in terms of the clusters of the eigenfrequencies in the entire
spectral band.

\vspace*{1cm}
\section{Introduction and statement of results}
\setcounter{equation}{0}
\setcounter{dref}{0}
In this paper we shall study some problems arising in the spectral
analysis of dissipative wave equations on compact manifolds. In 
order to formulate the results, we shall begin by recalling some standard assumptions and hypotheses. 

Let $M$ be a compact connected smooth Riemannian manifold of dimension
$\geq 2$, and let $\Delta$ be the corresponding 
Laplace-Beltrami operator. We consider the Cauchy problem for the wave equation with a damping term, 
\begeq
\label{1.1}
\left\{ \begin{array}{ll}
\left(-D_t^2-\Delta+2ia(x)D_t\right)u=0,\:\:(t,x)\in
\real_+\times M, \\
u |_{t=0}= u_0\in H^1(M),\quad D_t u|_{t=0}=u_1\in L^2(M).
\end{array} \right.
\endeq
Here $D_t=\partial_t/i$, and the damping coefficient $a$ is a bounded
nonnegative function on $M$, which is not identically zero. For
simplicity, we shall assume that $a\in C^{\infty}(M;\overline{\real_+})$. The motivation
for studying this problem comes from the geometric control theory (see~\cite{Lebeau}
and~\cite{AschLebeau}), where one is interested in the long time
behaviour of solutions to (\ref{1.1}), in relation to the
geometry of the underlying manifold and the control (damping) region.

Associated with the evolution problem (\ref{1.1}) is the solution
operator ${\cal U}(t)=e ^{it {\cal A}}$, $t\geq 0$, acting in the 
Hilbert space of the Cauchy data ${\cal H}=H^1\times L^2$ and mapping
$(u_0,u_1)\in {\cal H}$ to $(u(t,\cdot),D_tu(t,\cdot))$. Here we equip
${\cal H}$ with the norm 
$$
\left |\left|\,\left(\begin{array}{cc}
u_0 \\
u_1
\end{array}\right)\,\right |\right|_{{\cal H}}^2=\norm{\nabla
u_0}_{L^2}^2+\norm{u_0}_{L^2}^2+\norm{u_1}_{L^2}^2.
$$
The semigroup ${\cal U}(t)$
can be introduced by means of the Hille-Yosida theorem --- see Section 
2 for the corresponding resolvent estimates.  
The infinitesimal generator ${\cal A}$ is the operator
$$
{\cal A} = \left( \begin{array}{cc}
   0 & 1 \\
  -\Delta & 2ia(x)
\end{array} \right): {\cal H}\rightarrow {\cal H},
$$
with the domain $ D({\cal A})=H^2\times H^1$. Here
$H^s=H^s(M)$ is the standard Sobolev space on $M$. It follows that the
spectrum of ${\cal A}$, ${\rm Spec}({\cal A})$, is discrete. 

The energy of the solution $u(x,t)$ of (\ref{1.1}) at time $t$,  
\begeq
\label{1.2}
E(u,t)=\frac{1}{2}\int_M \left(\abs{\nabla u}^2+\abs{D_t u}^2\right)\,dx
\endeq
is nonincreasing as $t\rightarrow \infty$, and relations between the rate of the decay of the energy and the spectrum of 
${\cal A}$ were studied 
in~\cite{Lebeau} and~\cite{RT}; below we shall recall some of these
results. The asymptotic distribution of the eigenvalues 
of ${\cal A}$ has been studied by Sj{\"o}strand~\cite{Sjostrand1}, and some further results
in the case when the geodesic flow on $M$ is periodic have been
obtained in~\cite{Hitrik1}. We also remark that under the assumption
of periodicity of the flow, much more precise results are possible in the two-dimensional analytic case, and
we refer to~\cite{HSj1}, as well as to forthcoming papers in this series, for a
detailed study of this problem. In the present work we shall be
concerned with results, valid in all dimensions $\geq 2$.

We notice that $\tau\in \comp$ is an eigenvalue of ${\cal A}$ precisely when the equation
\begeq
\label{1.3}
\left(-\Delta+2i\tau a(x)-\tau ^2\right)u(x)=0
\endeq
has a nonvanishing smooth solution. The eigenvalues $\tau$ will also be called the eigenfrequencies. The multiplicity 
$m(\tau_0)\in \{1,2\ldots\}$ of an eigenfrequency $\tau_0$ is defined as the rank of the spectral projection
$$
\Pi_{\tau_0}=\frac{1}{2\pi i}\int_{\gamma} (\tau-{\cal A})^{-1}\,d\tau,
$$
where $\gamma$ is a sufficiently small circle centered at $\tau_0$. We
shall let $E_{\tau_0}:=\Pi_{\tau_0}{\cal H}$ stand for the generalized eigenspace corresponding to $\tau_0$. 

\smallskip
It follows easily from (\ref{1.3}) that if $\tau$ is an eigenfrequency, then
$$
0\leq \Im \tau\leq 2\norm{a}_{\infty},
$$
and since (\ref{1.3}) is invariant under the map $(\tau,u)\rightarrow
(-\overline{\tau},\overline{u})$, the set of the eigenfrequencies is
symmetric with respect to the reflection in the imaginary axis. We
remark that $\tau=0$ is an eigenfrequency corresponding to the constant
solution of (\ref{1.3}), and in~\cite{Lebeau} it was proved that
$\tau=0$ is the only eigenfrequency with the vanishing imaginary
part.    

Let us introduce a strongly continuous family of operators
\begeq
\label{1.3.1}
U(t): L^2\rightarrow H^1,\quad t\geq 0, 
\endeq
which takes a vector $f\in L^2$ to the first component of ${\cal U}(t)$ applied to a
vector of the form $(0,f)$. The mapping in (\ref{1.3.1})
will be referred to as the propagator. A simple computation shows that
in terms of $U(t)$, the matrix representation of the semigroup ${\cal U}(t)$ is given by 
$$
{\cal U}(t)=\left( \begin{array}{cc}
   V(t) & U(t) \\ 
   D_t V(t) & D_t U(t) \end{array} \right): {\cal H}\rightarrow {\cal H}, 
$$
where $V(t)f=D_t U(t)f-2i U(t)(af)$. 

In the self-adjoint case, i.e.\ when $a\equiv 0$, it is true that $U(t)=i\sin{t\sqrt{-\Delta}}/{\sqrt{-\Delta}}$, and we 
have a Fourier expansion,
$$
\frac{\sin{t\sqrt{-\Delta}}}{\sqrt{-\Delta}}f(x)=\sum_{\lambda_j^2\in \mathrm{Spec}(-\Delta)} 
\frac{\sin{\lambda_j t}}{\lambda_j} \varphi_j(x),\;\: -\Delta \varphi_j(x)=\lambda_j^2 \varphi_j(x),
$$
where the convergence is absolute in the case of smooth data. In this
paper, we shall be interested in eigenfunction expansions of 
the propagator $U(t)$, as $t\rightarrow \infty$, in the case of the non\-vani\-shing damping
term. We are therefore concerned with relations between the 
long time behaviour of ${\cal U}(t)$ and the spectrum of ${\cal
A}$. This problem is closely related to that of completeness and
summability of generalized eigenfunction expansions of the operator 
${\cal A}$. The classical theory of non-selfadjoint operators, as
described in~\cite{GoKr} and~\cite{Markus}, shows that the space of
finite linear combinations of the generalized eigenfunctions of ${\cal A}$, 
$\bigoplus_{\tau\in {\rm Spec}({\cal A})}E_{\tau}$, is dense in ${\cal H}$. However, 
as emphasized in [1], the lack of orthogonality in the direct sum
above indicates that the generalized eigenfunctions may still fail to form a (Riesz) basis in ${\cal H}$. 
The question of basisness in the context
of weak abstract non-selfadjoint perturbations of selfadjoint 
operators has been studied in~\cite{Markus}, and when specialized to
the case at hand, the results of~\cite{Markus} apply in dimension one
only, implying then the existence of a Riesz basis of the generalized
eigenfunctions. See also~\cite{CoxZuazua} and~\cite{Shubova}. 
Here we shall be concerned with the case $n\geq 2$, and as well as
using the classical non-selfadjoint theory, we shall also rely on
direct microlocal methods. 

In addition to~\cite{Markus}, in this paper we have been motivated by recent studies
of expansions for propagators on unbounded domains in terms of 
the resonance states---see~\cite{Zworski} for an overview of this work
and further references. Our Theorems 
1.2 and 1.3 below complement the corresponding results in the theory of
resonances in a situation, when the decay of states is caused by a direct 
dissipative mechanism on a compact manifold. In Theorem 1.4 we state a
result which, to the best of our 
knowledge, has no analogue for problems on unbounded domains. 

\vskip 2mm
When $p(x,\xi)=\xi^2$ is the principal symbol of $-\Delta$ defined on
$T^*M$, we introduce the Hamilton vector field $H_p$, and recall that the 
integral curves of the Hamilton flow $\exp(tH_p): p^{-1}(1)\rightarrow p^{-1}(1)$ are mapped via the projection 
$\pi: T^*M\rightarrow M$ to 
geodesics on $M$. When $T>0$, we put
\begeq
\label{1.3.5}
\hk{a}_T=\frac{1}{T} \int_0^T a(\exp(tH_p)\,dt\;\;\wrtext{on}\;\;p^{-1}(1). 
\endeq
We define also 
$$
0\leq A(T)=\inf_{p^{-1}(1)}\hk{a}_T, 
$$
which is a bounded continuous function on $\real_+$. We recall
from~\cite{Lebeau} that the limit 
\begeq
\label{1.4}
A(\infty):=\lim_{T\rightarrow \infty} A(T)
\endeq
exists and it is true that $A(\infty)=\sup_{T>0} A(T)$.  

\vskip 2mm
\noindent
In what follows we shall say that the geometric control condition 
holds if there exists a time $T_0>0$ such that any geodesic of length $\geq
T_0$ meets the open set $\{x; a(x)>0\}$. In this case it has been established by J.~Rauch and M.~Taylor~\cite{RT} that 
there is a uniform 
exponential decay of the energy $E(u,t)$. A proof of this statement can
also be extracted from the arguments given in Section 3 below. 

Our starting point is the following result, which is essentially due to~\cite{Lebeau}, and was also proved 
in~\cite{Sjostrand1}. 

\begin{theo} Assume that the geometric control condition holds. Then $A(\infty)>0$, and, for every 
$\varepsilon \in (0, A(\infty))$, there are at most finitely many
eigenfrequencies in the strip 
\begeq
\label{1.4.1}
\real+i[0,A(\infty)-\eps]. 
\endeq
\end{theo}

Our first result gives an expansion of the propagator $U(t)$, as
$t\rightarrow \infty$, in terms of the eigenfrequencies in the strip
(\ref{1.4.1}), provided that the geometric control condition is
satisfied. This results is analogous to the corresponding one in 
scattering theory, giving an expansion of the wave group in the
exterior of a nontrapping obstacle, in terms of exponentially decaying resonant
modes. See~\cite{TZ1} and references given there. 

\begin{theo}
Assume that the geometric control condition holds. Let $\varepsilon\in (0,A(\infty))$ be such that there are no 
eigenfrequencies with the imaginary part $=A(\infty)-\varepsilon$. Then there exists $C=C(\varepsilon)>0$ such that 
$$
U(t) = \sum_{\tau \in {\rm {Spec}}({\cal A}),\, {\rm Im} \tau
<A(\infty)-\varepsilon} e^{it\tau}p_{\tau}(t)+E(t), \quad t>0, 
$$
where $p_{\tau}(t)$ is a polynomial in $t$ with values in ${\cal L}(L^2,H^s)$
for any $s\in \real$, and the remainder term satisfies 
\begeq
\label{1.4.5}
\norm{E(t)}_{{\cal L}(L^2,L^2)}\leq C e^{-(A(\infty)-\varepsilon)t}.  
\endeq
We have 
$$
p_{\tau=0}(t)f=\frac{i\int f(x)\,dx}{\int 2a(x)\,dx},\quad f\in L^2. 
$$
\end{theo}

In order to motivate the next result, let us recall that
in~\cite{Lebeau}, G. Lebeau determined the optimal rate of the
exponential decay of the energy $E(u,t)$. See also~\cite{RT}. Set 
$$
D=\inf\{\Im \tau; \tau\in \Spec({\cal A}),\;\abs{\tau}>0\},
$$
and 
$$
\alpha=\sup\{\beta\geq 0; \;\exists B>0\;\forall u\in {\cal H},
\forall t\geq 0, E(u,t)\leq B e^{-\beta t} E(u,0)\}.
$$
Here the space ${\cal H}$ of the Cauchy data is identified with solutions
of (\ref{1.1}). It was proved in~\cite{Lebeau} that 
\begeq
\label{1.5}
\alpha=2\min(D,A(\infty)).
\endeq
In particular, it follows from (\ref{1.5}) that if there exists a
geodesic not intersecting the support of $a$, then
$\alpha=A(\infty)=0$. In~\cite{Lebeau}, an explicit example  
of a surface of revolution was given, where there exists a 
closed geodesic not meeting $\supp(a)$, while $D>0$, so that
the eigenfrequencies are bounded away from the real axis. In the 
example the geodesic in question is hyperbolic. Our next result
shows that if there exists a closed geodesic away from the damping
region, which is elliptic, then we have a sequence of eigenfrequencies
approaching rapidly the real axis. 

The precise assumptions on the geodesic will be formulated in
terms of the associated linearized Poincar\'e map, and we 
shall pause here to recall this notion. Given a simple closed $H_p$-trajectory
$\gamma: [0,T]\rightarrow p^{-1}(1)$ (satisfying
$\gamma(0)=\gamma(T)$), the linearized Poincar\'e map $P_{\gamma}$ is
defined as follows: let $H\subset p^{-1}(1)$ be a smooth hypersurface
intersecting $\gamma$ transversally at $\gamma(0)$. Then 
$P_{\gamma}$ is the differential at $\gamma(0)$ of the smooth locally
defined Poincare map: $H\rightarrow H$, obtained by following the flow of $H_p$ once around $\gamma$. 
Up to a symplectic conjugation, the mapping $P_{\gamma}$ can be viewed as a linear
symplectic transformation on $\real^{2(n-1)}$, and we say that $\gamma$ is
elliptic and nondegenerate if the spectrum 
of $P_{\gamma}$ lies on the unit circle and $1$ is not an eigenvalue. 

Given a nondegenerate elliptic closed $H_p$-trajectory $\gamma$, with
the eigenvalues of the corresponding $P_{\gamma}$ of the form $e ^{\pm i\alpha_j}$,
$j=1,\ldots n-1$, and $0<\alpha_j\leq \pi$, we say that $P_{\gamma}$ is $N$-ele\-men\-ta\-ry if 
$$
\sum_{j=1}^{n-1} k_j \alpha_j\neq 0\quad \mod\, 2\pi \z
$$
for each $k=(k_1,\ldots k_{n-1})\in \z^{n-1}$ with $0<\sum_{j=1}^{n-1}
\abs{k_j}\leq N$. 

\begin{theo}
Suppose that there exists a closed primitive $H_p$-orbit 
$\gamma: [0,T]\rightarrow p^{-1}(1)$, which is non\-de\-ge\-nerate 
el\-lip\-tic, and such that the ass\-ocia\-ted map\-ping $P_{\gamma}$ is
N-ele\-men\-tary for some $N\geq 4$. Assume that $\pi(\gamma([0,T]))\cap
\supp(a)=\emptyset$. Then there exists a sequence 
$(\tau_j)$, $j=1,2\ldots\,$, of eigenfrequencies such that $\abs{\Re \tau_j}\rightarrow \infty$, and
$0<\Im \tau_j\leq {\cal O}_N(1)\abs{\Re \tau_j}^{-N}$, $j\rightarrow \infty$, for any $N$. 
\end{theo} 

\Remark. This theorem can be viewed as an analogue of the results of
~\cite{StefanovVodev}, ~\cite{TZ}, and~\cite{Stefanov} in the theory
of resonances. 

\vskip 2mm
Our last theorem is concerned with an expansion of the propagator $U(t)$
on a Zoll manifold, i.e., a manifold all of whose geodesics are
closed. We refer to~\cite{Besse} for the background and motivation for the
study of such manifolds. In this case, the set of the eigenfrequencies
displays a cluster structure, and the expansion holds in terms of the disjoint clusters in the entire spectral band.

\begin{theo}
Let M be a Zoll manifold, i.e. assume that $\exp(\pi
H_p)(x,\xi)=(x,\xi)$, $(x,\xi)\in p^{-1}(1)$. Then
\begin{enumerate}
\item There exists $C>0$ such that all eigenfrequencies $\tau$, with
  $\Re\tau >0$, except for finitely many, are contained in the union
  of the rectangles 
\begeq
\label{1.6}
I_k=\left[k+\frac{\alpha}{4}-\frac{C}{k},k+\frac{\alpha}{4}+\frac{C}{k}\right]+
i[0,{\cal O}(1)], \quad k=1,2,\ldots,
\endeq
where $\alpha\in \z$ is the Maslov index of the closed
$H_p$-trajectory
$$
\{\exp(tH_p)(x,\xi);\,(x,\xi)\in p^{-1}(1),\,t\in [0,\pi]\}.
$$ 

\item We have 
\begeq
\label{1.7}
U(t) = \sum_{k=1}^{\infty} \left 
(\sum_{\tau \in {\rm Spec}({\cal A}),\,\tau \in I_k} e^{it\tau}
p_{\tau,k}(t)+
\sum_{\tau \in {\rm Spec}({\cal A}),\, -\overline{\tau} \in I_k} e^{it \tau} p_{\tau,k}(t)\right),\quad t>0.
\endeq
Here $p_{\tau,k}$ is a polynomial in $t$ with values in ${\cal
L}(L^2,H^s)$ for any $s\in \real$. The outer sum in {\rm (\ref{1.7})} converges absolutely in ${\cal L}(H^{\theta},L^2)$ for
every $\theta>0$ and for each fixed $t>0$. 
\end{enumerate}
\end{theo}

\Remark. The first part of Theorem 1.4 has been established in~\cite{Hitrik1}.
This is a well-known and in general optimal result in the self-adjoint case, 
$a\equiv 0$---see Chapter 29 of~\cite{HormIV}.
 
\Remark. The rectangles $I_k$ in (\ref{1.6}) become disjoint only for $k$ large enough,
$k\geq k_0\gg 1$, and the expansion in (\ref{1.7}) should be
understood in the following sense: If we set
$$
U_k(t)=\sum_{\tau \in {\rm Spec}({\cal A}),\,\tau \in I_k} e^{it\tau}
p_{\tau,k}(t) +\sum_{\tau\in {\rm Spec}({\cal A}),\,
-\overline{\tau}\in I_k} e^{it\tau} p_{\tau,k}(t), \quad k\geq k_0,
$$
then the series 
$$
\sum_{k=k_0}^{\infty} U_k(t)
$$
converges absolutely in ${\cal L}(H^{\theta},L^2)$, for any $\theta>0$ and each
$t>0$, and we have
$$
U(t)=\sum_{\tau \in {\rm Spec}({\cal A}),\,\abs{\Re \tau}<k_0+\alpha/4}
e^{it\tau}p_{\tau}(t)+\sum_{k \geq k_0} U_k(t) \quad
\wrtext{in}\,\,{\cal L}(H^{\theta},L^2).
$$

\vskip 6mm
\noindent
The plan of the paper is as follows. After some preliminaries in
Section 2, Theorem 1.2 is proved in Section 3. 
The proof relies on a simple conjugation by pseudodifferential operators, microlocally in the region of 
the phase space, 
where we lack ellipticity. We also adapt the classical argument of Morawetz~\cite{Mor} to our situation, 
in order to obtain Theorem 1.2 in the 
case of the $L^2$ initial data. In Section 4 we study the case when the geodesic flow is periodic. Here we recall the 
methods of~\cite{Hitrik1}, and 
show how they can be used to prove the second part of Theorem 1.4. Section 5 is devoted to the proof of Theorem 1.3. 
The proof depends on an a priori semiclassical resolvent estimate, and here we follow the idea of~\cite{Markus} 
and~\cite{Sjostrand1} of creating gaps in 
the spectrum of the unperturbed operator by means of an auxiliary finite
rank perturbation. At this stage we also use the results and methods
of~\cite{PopovCardoso},~\cite{TZ1}, and~\cite{Stefanov}. 

\vskip 2mm
The results of this paper have been announced in~\cite{Hitrik2}. 

\vskip 2mm
\noindent
{\it Acknowledgements.} I am very grateful to Maciej Zworski for numerous very helpful discussions on the subject of 
the present paper. I would also like to thank Nicolas Burq for
pointing out the argument of Morawetz, a modification 
of which is used in Section 3. The support of the Swedish Foundation
for International Cooperation in Research and Higher Education (STINT), as 
well as of the MSRI postdoctoral fellowship, is gratefully acknowledged.  

\section{Preliminaries}
\setcounter{equation}{0}
In the beginning of this section we shall recall, for future
reference, some resolvent estimates for ${\cal A}$. 
In doing so we shall write $D=\{\tau \in \comp; \Im \tau \in [0,
2\norm{a}_{\infty}]\}$, and we recall that ${\rm Spec}({\cal A})\subset
D$. 

\begin{prop}
There exists a constant $C>0$ such that
\begeq
\label{2.1}
\norm{\left(\tau-{\cal A}\right)^{-1}}_{{\cal L}({\cal H},{\cal H})}\leq \frac{C}
{\dist(\tau, D)},\quad \tau \not \in D. 
\endeq
\end{prop}
\begin{proof}
Let us consider the following closed codimension one subspace of ${\cal H}=H^1\times L^2$,
$$
{\cal H}_0=\{(u_0,u_1)\in {\cal H}; \int 2ia(x)u_0(x)\,dx-\int u_1(x)\,dx=0\}.
$$
Since $a\geq 0$ is not identically zero, it is true that the inner product on
${\cal H}_0$,
$$
\left<(u_0,u_1), (v_0,v_1)\right>_E:=\int\!\nabla u_0 \overline{\nabla v_0}\,dx+
\int\!u_1\overline{v_1}\,dx
$$
defines a norm $\norm{\cdot}_E$ there which is equivalent to the
natural norm on ${\cal H}$, $\norm{\cdot}_{{\cal H}}$. When $U=(u_0,u_1)\in D({\cal A})$, it 
is also true that 
${\cal A}U=(u_1, -\Delta u_0+2ia u_1)$ is in ${\cal H}_0$, and a simple computation shows that ${\cal H}_0$ contains 
$\bigoplus_{0\neq \tau \in {\rm Spec}({\cal A})} E_{\tau}$ as a (dense) subspace. When $U\in D({\cal A})\cap {\cal H}_0$, we have 
$$
0\leq \Im \left<{\cal A}U,U\right>_E\leq 2\norm{a}_{\infty}\norm{U}^2_E,
$$
so that when $\Im \tau<0$, we get $\Im \left<\left({\cal
A}-\tau\right)U,U\right>_E\geq \abs{\Im \tau} 
\left<U,U\right>_E$. 
It follows that there exists $C>0$ such that for any $U\in {\cal H}_0$, 
\begeq
\label{2.2}
\norm{(\tau-{\cal A})^{-1}U}_{{\cal H}}\leq
\frac{C}{\abs{\Im \tau}}\norm{U}_{{\cal H}},\quad \Im \tau<0.
\endeq
Now it is easily seen that the space $E_0$ is of dimension one and is generated by $(1,0)$. We have 
${\cal H}={\cal H}_0\oplus E_0$, 
also in the topological sense, and since 
\begeq
\label{2.3}
\norm{(\tau-{\cal A})^{-1}V}_{{\cal H}}\leq
\frac{1}{\abs{\Im \tau}}\norm{V}_{{\cal H}},\quad
\Im \tau<0,\;\;V\in E_0,
\endeq
we conclude that (\ref{2.1}) follows from (\ref{2.2}) and (\ref{2.3}), in the case when $\Im \tau <0$. 
The proof is complete, for in the case 
when $\Im \tau> 2\norm{a}_{\infty}$, we argue in a similar way.
\end{proof}

\Remark. A straightforward computation shows that 
\begeq
\label{2.4}
(\tau-{\cal A})^{-1}=\left( \begin{array}{cc}
   R(\tau)(2ia-\tau) & -R(\tau) \\
   R(\tau)(2ia\tau-\tau^2)-1 & -R(\tau)\tau
\end{array} \right),
\endeq 
where 
\begeq
\label{2.4.5}
R(\tau)=P(\tau)^{-1},\quad P(\tau)=-\Delta+2ia\tau-\tau^2.
\endeq 
Combining (\ref{2.4}) with Proposition 2.1 we get
\begeq
\label{2.5}
\norm{R(\tau)}_{{\cal L}(L^2,L^2)}\leq \frac{C}{\abs{\tau}\abs{\Im \tau}},\quad \Im \tau<0. 
\endeq
and 
\begeq
\label{2.5.5}
\norm{R(\tau)}_{{\cal L}(L^2,L^2)}\leq \frac{C}{\abs{\tau}(\Im
\tau-2\norm{a}_{\infty})},\quad \Im \tau> 2\norm{a}_{\infty}.
\endeq

\vskip 2mm
It follows from the proof of Proposition 2.1 that 
$$
\norm{\left(\tau-{\cal A}\right)^{-n}}_{{\cal L}({\cal H},{\cal H})}
\leq \frac{C}{\dist(\tau, D)^n},\quad \tau \not \in D,\,\, n=1,2,\ldots, 
$$
and an application of the Hille-Yosida theorem shows that ${\cal A}$ generates a uniformly bounded strongly continuous 
semigroup ${\cal U}(t)=e^{it {\cal A}}$, $t\geq 0$, in ${\cal H}$, so that in particular 
$U(t)\in C(\overline{\real_+}, {\cal L}(L^2,H^1))\cap C^1(\overline{\real_+},{\cal L}(L^2,L^2))$, in the strong sense. 
See~\cite{Davies}. We 
have the following relation between ${\cal U}(t)$ and the resolvent of ${\cal A}$, valid when $x\in D({\cal A}^2)$, 
\begeq
\label{2.6}
e^{it{\cal A}} x=x+\frac{1}{2\pi i}\int_{-\infty-i\alpha}^{\infty-i\alpha}
e^{it\tau} \left(\tau-{\cal A}\right)^{-1} {\cal A}x \,\frac{d\tau}{\tau},
\quad t>0,\;\alpha>0.
\endeq
The integral in (\ref{2.6}) converges absolutely in ${\cal H}$. If we
choose $x=(0,f)$, $f\in H^2$, then (\ref{2.4}) 
shows that the first component of $\left(\tau-{\cal A}\right)^{-1} {\cal A}x$ is equal to $-\tau R(\tau)f$. Thus, 
with the absolute convergence in $L^2$ we have,   
\begeq
\label{2.7}
U(t)f =\frac{i}{2\pi}\int_{-\infty-i\alpha}^{\infty-i\alpha}
e^{it\tau} R(\tau) f \,d\tau,\quad \alpha>0,\;\;\quad f\in H^2(M). 
\endeq
Notice that when $f\in L^2(M)$, (\ref{2.7}) still holds, with the integral in (\ref{2.7}) 
converging in $L^2(\Im \tau=-\alpha; L^2(M))$. 

\vskip 4mm
When establishing our results, it will be convenient to work in the
semiclassical setting. We shall therefore perform a semiclassical
reduction, following~\cite{Sjostrand1} and~\cite{Hitrik1}. When considering the eigenvalue problem 
$$
(-\Delta+2ia(x)\tau-\tau^2)u(x)=0,
$$
with $\Im \tau={\cal O}(1)$ and $\Re \tau \gg 1$, we write $\tau=\sqrt{z}/h$, where $0<h\ll 1$,
and $z$ belongs to the fixed domain $[\alpha,\beta]+i[-\gamma,\gamma]$,
$0<\alpha<1<\beta<\infty$, $\gamma>0$. We are then led to study the
$h$-differential operator 
\begeq
\label{2.8}
{\cal P}=P+ihQ(z),
\endeq
where $P=-h^2 \Delta$ and $Q(z)=2a(x)\sqrt{z}$. We notice that $P$ is
essentially selfadjoint on $C^{\infty}(M)$, with the domain of the
selfadjoint realization being $H^2(M)$, and $Q(z)$ is bounded and
selfadjoint for $z$ real positive. 

In the main part of the paper, we shall make use of some elementary
semiclassical pseudodifferential calculus, and we digress here to
recall some relevant notation. Let $S^m(T^*M)=S^m_{1,0}(T^*M)$ be the
space of functions $a(x,\xi;h)$ on $T^*M\times(0,h_0]$, $h_0>0$, which
are smooth in $(x,\xi)$ and such that  
\begeq
\label{2.9}
\partial_x^{\alpha}\partial_{\xi}^{\beta}a(x,\xi;h)=
{\cal O}_{\alpha,\beta}(1)\left(1+\abs{\xi}\right)^{m-\abs{\beta}}, 
\quad (x,\xi)\in T^*M,\,\,h\in (0,h_0]. 
\endeq
When $a$ depends on some additional parameters, we require the estimates (\ref{2.9}) 
to hold uniformly with respect to these parameters. The formula for the classical $h$-quantization, 
$$
a(x,hD_x;h)u(x)=
\frac{1}{(2\pi h)^n}\int\!\!\int e^{i\frac{\hk{x-y,\xi}}{h}} a(x,\xi;h) u(y)\,dy\,d\xi,
$$
then defines a class of $h$-pseudodiffential operators on $M$, which
will be denoted by $\wrtext{Op}_h(S^m)$.

For $h$-dependent symbols, we write $a\in S_{\mathrm{cl}}^m$ if there exists $a_0\in S^m$, independent of 
$h$, such that $a-a_0\in h S^{m-1}$, and we say that $a_0$ is the principal symbol of $a$ and of
the corresponding operator. When $A\in \wrtext{Op}_h(S^m)$ has a principal symbol $a_0\in S^m$, 
we write $A=\wrtext{Op}_h(a_0)$ and say that $A$ is classical.

We shall finally recall the notion of the (semiclassical) wave front set, $WF_h(A)$, of an operator 
$A\in \wrtext{Op}_h(S^m)$. When $\rho \in T^*M$, we say that $\rho \not \in WF_h(A)$, if the full symbol 
of $A$, for some choice of local coordinates near the projection of $\rho$, is of class 
$S^{-\infty,m}:=\cap_{k} h^k S^m$ near $\rho$. It follows from the definition that $WF_h(A)$ is a closed 
subset of $T^*M$, $WF_h(AB)\subset WF_h(A)\cap WF_h(B)$, and
$WF_h(A)=\emptyset \Rightarrow h^{-k}A\in \wrtext{Op}_h(S^m)$, for any
$k$. In this paper, we shall work in a purely semiclassical setting,
and consequently, only compact subsets of $T^*M$ will be
important. The definition of $WF_h$ given above will therefore be
sufficient for our purposes. 

\section {Proof of Theorem 1.2}
\setcounter{equation}{0}

The main step in the proof of Theorem 1.2 consists of establishing the following result. 
 
\begin{prop}
Assume that the Rauch-Taylor condition holds. Then for every
$\varepsilon\in (0,A(\infty))$ there exists $C=C(\varepsilon)>0$ such that if $\Im\tau\in [-1, A(\infty)-\varepsilon]$ and 
$\abs{\Re \tau}\geq C$ then 
\begeq
\label{3.1}
\norm{R(\tau)}_{{\cal L}(L^2,L^2)}\leq \frac{C}{1+\abs{\tau}}.
\endeq
\end{prop}
\begin{proof}
Performing the semiclassical reduction of Section 2, we can
reformulate the statement of the proposition in the following form: 
for every $\varepsilon\in (0,A(\infty))$ there exists $C>0$ and $h_0>0$ such that for
all $h\in (0,h_0]$ and all 
\begeq
\label{3.4}
z\in \Omega(\eps):=\{\Re z\in [\alpha,\beta], -2\sqrt{\Re z}\leq \Im z/h 
\leq 2\sqrt{\Re z}\left(A(\infty)-\varepsilon\right)\},
\endeq 
it is true that 
\begeq
\label{3.5}
\norm{({\cal P}-z)^{-1}}_{{\cal L}(L^2,L^2)}\leq \frac{C}{h}.
\endeq
Here ${\cal P}$ has been defined in (\ref{2.8}). The proof of
(\ref{3.5}), sketched in~\cite{Hitrik2} was a contradiction argument, 
and it relied upon a propagation estimate for 
a suitable semiclassical measure. Here we shall give a different proof, which has the merit of being more 
direct. In doing so, we shall use some arguments of~\cite{Sjostrand1}
and~\cite{Hitrik1}. 

In what follows we shall assume, as we may, that $\Im z={\cal
O}(h)$. Following the argument of Section 2 of~\cite{Sjostrand1}, we conjugate
${\cal P}=P+ihQ(z)$ by an ell\-ip\-tic self\-adjo\-int
$h$-pseudo\-diffe\-ren\-tial ope\-ra\-tor $A=\wrtext{Op}_h(e^g)$,
where $g\in S^{-1}(T^*M)$ is to be chosen. We have 
\begeq
\label{3.6}
A^{-1}\left(P+ihQ(z)\right)A=P+ih \wrtext{Op}_h(q(\Re z)-H_p g)+h^2 R(z). 
\endeq
Here $q(z)=2a\sqrt{z}$ and $R(z)\in \wrtext{Op}_h(S^0)$. Let us choose some $\delta>0$, sufficiently small but 
fixed. Continuing to follow the 
argument of~\cite{Sjostrand1}, we see that we can choose $g=g(\Re z)\in S^{-1}$ such that 
$$
q(\Re z)-H_p g=\langle{q(\Re z)}\rangle_T \quad \wrtext{on}\,\, p^{-1}\left([\alpha-\delta,\beta+\delta]\right).
$$
Here the trajectory average $\hk{q(\Re z)}_T$ is given by 
$$
\hk{q(\Re z)}_T=\frac{1}{T}\int_0^T q(\Re z)\circ \exp(tH_p)\,dt,
$$
and $T>0$ is to be chosen later, large enough but fixed. We then consider the right hand side of (\ref{3.6}), 
\begeq
\label{3.7}
\widehat{{\cal P}}=\widehat{{\cal P}}_T=P+ih \wrtext{Op}_h(\widehat{q}_T)+h^2 R_T(z), 
\endeq
where $\widehat{q}_T\in S^0$ is such that 
$$
\widehat{q}_T=\langle{q(\Re z)}\rangle_T\quad \wrtext{on}\,\, p^{-1}([\alpha-\delta,\beta+\delta]). 
$$

We now claim that 
$$
\lim_{T\rightarrow \infty}\inf_{p^{-1}(E)}\langle{q(\Re z)}\rangle_T=2\sqrt{\Re z}A(\infty), 
$$
locally uniformly in $E>0$. Indeed, we observe that the homogeneity property of the $H_p$-flow implies that 
$$
\inf_{p^{-1}(E)}\langle{a}\rangle_T=\inf_{p^{-1}(1)}\langle{a\rangle}_{\sqrt{E}T},
$$
and choosing $T$ large enough but fixed, depending on $\eps>0$, we can therefore
arrange that for $\widehat{q}_T$ in (\ref{3.7}), we have 
$$
\widehat{q}_T\geq 2\sqrt{\Re z}\left(A(\infty)-\eps/2\right)\quad \wrtext{on}\,\, p^{-1}([\alpha-\delta, \beta+\delta]). 
$$
Take now a cutoff function $0\leq\chi\in C^{\infty}_0(p^{-1}((\alpha-\delta,\beta+\delta)))$, such that $\chi=1$ near 
$p^{-1}([\alpha,\beta])$. An application of the G\aa{}rding inequality
allows us to conclude that for 
$z\in \Omega(\eps)$, 
\begeq
\label{3.8}
\left(\Im (\widehat{\cal P}-z)\chi u|\chi u\right)\geq \eta h\norm{\chi u}^2-{\cal O}(h^{\infty})\norm{u}^2,
\endeq
where $\eta=\eta(\eps)>0$ is fixed and $h>0$ is small enough, depending on
$\eps$. Here we let $\chi$ also stand for the corresponding quantization, and the inner
product and the norm in (\ref{3.8}) are taken in $L^2(M)$. On the other hand, we have 
$$
\left(\Im (\widehat{{\cal P}}-z)\chi u|\chi u\right)=\Im \left((\chi(\widehat{{\cal P}}-z)u|\chi u)+
([\widehat{{\cal P}},\chi]u|\chi u)\right),
$$
and the absolute value of this expression does not exceed
\begeq
\label{3.9}
{\cal O}(1) \norm{(\widehat{{\cal P}}-z)u}\,\norm{\chi u}+{\cal O}(1)\norm{[\widehat{{\cal P}},\chi]u}\,\norm{\chi u}. 
\endeq
Let now $\chi_1\in C^{\infty}_0(p^{-1}((\alpha-\delta,\beta+\delta))$ be such that $\chi_1=1$ near 
$p^{-1}([-\alpha,\beta])$ and 
the support of $\chi_1$ is contained in the interior of the set where $\chi=1$. Then 
$$
\norm{[\widehat{{\cal P}},\chi]u}\leq \norm{[\widehat{{\cal P}},\chi](1-\chi_1)u}+{\cal O}(h^{\infty})\norm{u} 
\leq {\cal O}(h)\norm{(1-\chi_1)u}+{\cal O}(h^{\infty})\norm{u}. 
$$
Now the operator $\widehat{\cal P}-z$ is elliptic away from $p^{-1}([\alpha,\beta])$, and the semiclassical elliptic 
regularity shows that 
\begeq
\label{3.10}
\norm{(1-\chi_1)u}\leq {\cal O}(1)\norm{(\widehat{\cal P}-z)u}+{\cal O}(h^{\infty})\norm{u}.
\endeq
We get
\begeq
\label{3.11}
\norm{[\widehat{{\cal P}},\chi]u}\leq {\cal O}(h)\norm{(\widehat{\cal P}-z)u}+{\cal O}(h^{\infty})\norm{u}.
\endeq
Combining (\ref{3.8}) together with (\ref{3.9}) and (\ref{3.11}), we infer that 
\begin{eqnarray*}
\norm{\chi u}^2 & \leq & \frac{{\cal O}(1)}{h}\norm{(\widehat{\cal
P}-z)u}\,\norm{\chi u}+
{\cal O}(h^{\infty})\norm{u}^2 \\
& \leq & \frac{{\cal O}(1)}{h^2}\norm{(\widehat{\cal P}-z)u}^2
+\frac{1}{2}\norm{\chi u}^2+
{\cal O}(h^{\infty})\norm{u}^2,
\end{eqnarray*}
and therefore
$$
\norm{\chi u}\leq \frac{{\cal O}(1)}{h}\norm{(\widehat{\cal P}-z)u}+{\cal O}(h^{\infty})\norm{u}.
$$
Combining this with (\ref{3.10}), with $\chi$ in place of $\chi_1$, we get
$$
\norm{u}\leq \frac{{\cal O}(1)}{h}\norm{(\widehat{\cal P}-z)u},\quad u\in H^2(M),\quad z\in \Omega(\eps). 
$$
This completes the proof if we take into account the uniform bounds
$$
A,\,A^{-1}={\cal O}_{s}(1): H^s \rightarrow H^s,\quad s\in \real,  
$$
where the $H^s$-spaces have been equipped with their natural
$h$-dependent norms. 
\end{proof}

\smallskip
It is now easy to obtain an expansion of $U(t)$ in ${\cal
L}(H^{\theta},L^2)$, $\theta>0$, as $t\rightarrow \infty$. 
Indeed, if $f\in H^2$ then 
$$
R(\tau)(-\Delta+2ia\tau-\tau^2)f=f,
$$
which together with Proposition 3.1 implies that for $\abs{\Re \tau}\geq C$ and $\Im \tau\in [-1,A(\infty)-\varepsilon]$, we have 
\begeq
\label{3.10.5}
\norm{R(\tau)}_{{\cal L}(H^2,L^2)}\leq \frac{{\cal O}(1)}{(1+\abs{\tau})^2}.
\endeq
Interpolating between (\ref{3.1}) and (\ref{3.10.5}), we get  
$$
\norm{R(\tau)}_{{\cal L}(H^{\theta},L^2)}\leq \frac{{\cal O}(1)}{(1+\abs{\tau})^{1+{\theta}/2}}, \quad 0\leq \theta \leq 2, 
$$
so that for $\theta >0$, 
$$
\int_{-\infty-i\alpha}^{\infty-i\alpha}e^{it\tau} R(\tau)\,d\tau \in
C(\overline{\real}_+,{\cal L}(H^{\theta},L^2)),\quad \alpha>0.
$$
The expansion for $U(t)$ stated in Theorem 1.2, with an error estimate
(\ref{1.4.5}) in ${\cal L}(H^{\theta},L^2)$, 
is then obtained by deforming the contour of integration in the right hand side of 
(\ref{2.7}), exactly as in the proof of Theorem 1 of~\cite{TZ1}. 

\vskip 2mm
\noindent
We now come to prove Theorem 1.2 in the Banach space ${\cal
L}(L^2,L^2)$. In doing so, we shall make use of an argument 
which comes essentially from~\cite{Mor}, there given in the context of the wave equation in the exterior 
of a nontrapping obstacle. See also~\cite{Burq}. We shall transform the initial value problem (\ref{1.1}) into a 
non-homogeneous problem.

Let $\eps\in (0,A(\infty))$ be such that there are no eigenfrequencies $\tau$ with the imaginary part 
$\Im \tau=A(\infty)-\eps$, and let 
$\tau_1,\ldots \tau_N$ be the eigenfrequencies with $\Im \tau_j < A(\infty)-\eps$. We introduce the 
spectral projector onto $\bigoplus_{j=1}^N E_{\tau_j}$, 
$$
\Pi =\frac{1}{2\pi i}\int_{\gamma} \left(\tau-{\cal A}\right)^{-1}\,d\tau,
$$
where $\gamma$ consists of the union of small positively oriented circles centered at the $\tau_j$'s. When 
$x\in {\cal H}$, it is true that 
$$
\left(\tau-{\cal A}\right)^{-1}\left(I-\Pi\right)x
$$
is a holomorphic function in a neigbourhood of the set where $\Im \tau\leq A(\infty)-\eps$, with values in 
${\cal H}$. We shall now derive weighted decay estimates for
$e^{it{\cal A}}(I-\Pi)x$. When doing so, we notice 
first that Proposition 3.1 together with an interpolation argument shows that 
\begeq
\label{3.12}
R(\tau)={\cal O}(1): L^2\rightarrow H^1, \quad \abs{\Re \tau}\geq {\cal O}(1)
\endeq
when $\Im \tau \leq A(\infty)-\eps$. Let $\psi \in C^{\infty}(\real;\real)$ be an increasing function such that 
$\psi(t)=0$, $t\leq 0$, and $\psi(t)=1$, $t\geq 1$. We consider the function $\psi(t)u(t,x)$ where $u(t,x)$ solves 
the evolution problem (\ref{1.1}) with the initial data $(u_0,u_1)=(I-\Pi)x$, when $x=(0,f)$, $f\in L^2$. Then 
\begeq
\label{3.13}
\left(-D_t^2-\Delta+2ia(x)D_t\right)\psi(t)u(t,x)=g(t,x), 
\endeq
where $g\in C(\real;L^2)$ is given by 
$$
g(t,x)=\psi'(t)\left(2a(x)u(t,x)+iD_tu(t,x)\right)+iD_t\left(\psi'(t)u(t,x)\right).
$$
The support of $g(\cdot, x)$ is contained in $[0,1]$ so that $g\in L^2(\real;L^2)$, and it follows from the 
Hille-Yosida theorem that 
\begeq
\label{3.13.1}
\norm{g}_{L^2({\rm {\bf R}};L^2)}\leq {\cal O}(1) \norm{f}_{L^2}. 
\endeq
Taking the Fourier transform in $t$ of (\ref{3.13}), ${\cal
F}_{t\rightarrow \tau}$, we get 
\begeq
\label{3.13.2}
P(\tau) {\cal F}_{t\rightarrow \tau} \psi(t)u(t,x)={\cal F}_{t\rightarrow \tau} g(t,x),
\quad {\cal F}_{t\rightarrow \tau}g(t,x) \in L^2(\real;L^2).
\endeq
Now an application of (\ref{2.4}) shows that 
$R(\tau) {\cal F}_{t\rightarrow \tau}g(t,x)=P(\tau)^{-1}{\cal
F}_{t\rightarrow \tau}g(t,x)$ is the first component of 
$$
\frac{1}{i}\left(\tau-{\cal A}\right)^{-1} {\cal F}_{t\rightarrow \tau} \psi'(t) \left( \begin{array}{cc}
   u(t,x) \\
   D_t u(t,x)
\end{array} \right),
$$
so that it depends holomorphically on $\tau$ in the set where $\Im \tau\leq
A(\infty)-\eps$, and takes values in $H^1$. 

An application of Parseval's formula to (\ref{3.13.2}) together with (\ref{3.12}) and
(\ref{3.13.1}) gives that 
\begin{eqnarray}
\label{3.14}
& & \norm{e^{(A(\infty)-\eps)t}\psi(t)u(t)}_{L^2({\rm {\bf R}};
H^1)} \\ \nonumber 
& = & \norm{R(\tau) {\cal F}_{t\rightarrow \tau}
g(t,\cdot)|_{\tau=\cdot+i(A(\infty)-\eps)}}_
{L^2({\rm {\bf R}}; H^1)} \\ \nonumber
& \leq & {\cal O}(1) 
\norm{{\cal F}_{t\rightarrow \tau}
g(t,\cdot)|_{\tau=\cdot+i(A(\infty)-\eps)}}_{L^2({\rm {\bf
R}};L^2)}\leq {\cal O}(1)\norm{f}_{L^2}. 
\end{eqnarray}
This gives an exponential decay estimate in the integrated form. When deri\-ving estimates that are pointwise in time, 
we consider the non-homogeneous equation 
satisfied by $\psi_T(t)u(t,x)$, where $\psi_T(t)=\psi(t-T)$, $T>1$. We
have
$$
\left(-D_t^2-\Delta+2ia(x)D_t\right) \psi_T(t)u(t,x)=g_T(t,x),
$$
where
$g_T(t,x)=(\psi_T''(t)+2a(x)\psi_T'(t))u(t,x)+2i\psi_T'(t)D_t u(t,x)$. Standard hyperbolic estimates (see~\cite{HormIV}) show that 
$$
\norm{u(T+1,\cdot)}_{L^2} \leq {\cal O}(1) \int_T^{T+1}\norm{g_T(t,\cdot)}_{L^2}\,dt \leq 
{\cal O}(1) \left(\int_T^{T+1}\norm{g_T(t,\cdot)}_{L^2}^2\,dt\right)^{1/2}, 
$$
and therefore
\begeq
\label{3.15}
\norm{u(T+1,\cdot)}_{L^2}^2 \leq {\cal O}(1) \left(\int_T^{T+1} \left(\norm{u(t,\cdot)}_{L^2}^2 +
\norm{D_t u(t,\cdot)}_{L^2}^2 \right)\,dt\right). 
\endeq
Now it follows from (\ref{3.14}) that 
\begin{eqnarray}
\label{3.16}
\int_T^{T+1} \norm{u(t,\cdot)}^2_{L^2}\,dt & \leq & e^{-2(A(\infty)-\eps)T}\int_0^{\infty} 
e^{2(A(\infty)-\eps)t}\psi^2(t)\norm{u(t,\cdot)}^2_{L^2}\,dt \\ \nonumber
& \leq & {\cal O}(1) e^{-2(A(\infty)-\eps)T}\norm{f}^2_{L^2},
\end{eqnarray}
and we shall now verify that the corresponding estimate holds when $u(t,\cdot)$ is replaced by $D_t u(t,\cdot)$. 
In doing so, we observe that (\ref{3.13.2}) gives
$$
P(\tau) {\cal F}_{t\rightarrow \tau} \left(D_t(\psi(t)u(t,x))\right)=\tau
{\cal F}_{t\rightarrow \tau} g(t,x),
$$
and since
$$
\tau R(\tau)={\cal O}(1): L^2 \rightarrow L^2, \quad \abs{\Re \tau}\geq {\cal O}(1), 
\Im \tau \in [-1, A(\infty)-\eps],
$$
repeating the previous arguments we see that 
$$
\int_0^{\infty} e^{2(A(\infty)-\eps)t}\norm{D_t(\psi(t)u(t,\cdot))}_{L^2}^2 \,dt \leq {\cal O}(1) \norm{f}^2_{L^2}.
$$
Since $\psi'$ is supported in $[0,1]$, we get 
$$
\int_0^{\infty} e^{2(A(\infty)-\eps)t}\psi^2(t)\norm{D_t u(t,\cdot)}_{L^2}^2 \,dt \leq {\cal O}(1) \norm{f}^2_{L^2},
$$
and therefore 
\begeq
\label{3.17}
e^{2(A(\infty)-\eps)T}\int_T^{T+1} \norm{D_t u(t,\cdot)}^2_{L^2}\,dt \leq {\cal O}(1) \norm{f}^2_{L^2}.
\endeq
Combining (\ref{3.15}) with (\ref{3.16}) and (\ref{3.17}) we conclude that 
\begeq
\label{3.18}
\norm{u(t,\cdot)}_{L^2} \leq {\cal O}(1) e^{-(A(\infty)-\eps)t}\norm{f}_{L^2}. 
\endeq

In order to derive an expansion of Theorem 1.2, we now only have to consider 
$$
e^{it{\cal A}}\Pi x=\frac{1}{2\pi i}\int_{-\infty-i\alpha}^{\infty-i\alpha}
e^{it\tau}\left(\tau-{\cal A}\right)^{-1}\Pi x\,d\tau, \quad
\alpha>0,\quad x=(0,f), 
$$
which is given by the sum of the residues of the integrand at the poles $\tau_j$, $j=1,\ldots N$. If 
$\tau_0$ is such a pole, then near $\tau_0$, 
we can write for some $M \in \nat$, not exceeding the multiplicity of
$\tau_0$, 
$$
\left(\tau-{\cal A}\right)^{-1} \equiv \frac{\Pi_0}{\tau-\tau_0}+\sum_{j=1}^M\frac{D_j}{(\tau-\tau_0)^{j+1}}, 
$$ 
modulo a function holomorphic in a \neigh{} of $\tau_0$. Here $\Pi_0$ is the spectral projector onto 
$E_{\tau_0}$ and 
$D_j=({\cal A}-\tau_0)^j \Pi_0$. The corresponding contribution to the
sum of the residues is of the form $e^{it\tau_o}$ times a polynomial of
degree at most $M$ in $t$ with values in 
${\cal L}({\cal H}, D({\cal A}^{\infty}))$. In particular, when $\tau_0=0$, it is
true that the corresponding generalized eigenspace $E_{0}$ is 
spanned by $(1,0)$ so that the pole at $0$ is simple, and the spectral
projector $\Pi_0$ is given by 
$$
\Pi_0 \left( \begin{array}{ccc}
u_0 \\
u_1
\end{array} \right)=\left( \begin{array}{ccc}
c \\
0
\end{array} \right),
$$
where
$$
c=\frac{\int 2ia(x)u_0(x)\,dx-\int u_1(x)\,dx}{\int 2ia(x)\,dx}.
$$
Combining these remarks together with (\ref{3.18}), we obtain the
statement of Theorem 1.2. 

\section{Proof of Theorem 1.4}
\setcounter{equation}{0}

In this section we shall work under the assumption that the geodesic flow on $M$ is periodic. To fix the ideas,
we shall assume throughout that $\exp(\pi H_p)(x,\xi)=(x,\xi)$,
$(x,\xi)\in p^{-1}(1)$.  
Our purpose is to derive an expansion for the damped wave propagator
$U(t)$ on $M$. In doing so, 
we shall make use of the methods of~\cite{Hitrik1}. Let us recall that
in~\cite{Hitrik1} it was established that 
$$
\Spec({\cal A})\cap \{\tau\in \comp; \Re \tau>0\}\subset \bigcup_{k=1}^{\infty} I_k,
$$
where
$$
I_k=\left[k+\frac{\alpha}{4}-\frac{{\cal
      O}(1)}{k},k+\frac{\alpha}{4}+\frac{{\cal
      O}(1)}{k}\right]+i[0,{\cal O}(1)], \quad \alpha\in \z, \quad k=1,2,\ldots.
$$

We shall now complement this result by deriving resolvent estimates in the gaps between the rectangles $I_k$. 

\begin{prop}
There exists $C>0$ such that when $\Im \tau\in [-1, 2\norm{a}_{L^{\infty}}+1]$ we have 
\begeq
\label{4.1}
\norm{R(\tau)}_{{\cal L}(L^2,L^2)}\leq \frac{C}{1+\abs{\tau}}, \quad
\abs{\Re \tau}\geq C,\;\;\dist\left(\abs{\Re \tau}, \cup_{k=1}^{\infty} J_k)\right)\geq \frac{1}{C}.
\endeq
Here $J_k=\left[k+\frac{\alpha}{4}-\frac{C}{k},k+\frac{\alpha}{4}+\frac{C}{k}\right]$.
\end{prop}

\begin{proof}
Applying the semiclassical reduction of Section 2, we see that in
order to prove (\ref{4.1}) it suffices to establish the existence 
of $C>0$ and $h_0>0$ such that for all $h\in (0,h_0]$ we have that 
\begeq
\label{4.2}
\norm{({\cal P}-z)^{-1}}_{{\cal L}(L^2,L^2)}\leq \frac{C}{h},
\endeq
when $\Im z={\cal O}(h)$ and $\dist\left(\Re z, \cup_{k=1}^{\infty} J_k(h)\right)\geq h/C$. Here the intervals
$$
J_k(h)=\left[h^2\left(k+\frac{\alpha}{4}\right)^2-{\cal O}(h^2),
h^2\left(k+\frac{\alpha}{4}\right)^2+{\cal O}(h^2)\right],\quad k\sim
\frac{1}{h}, 
$$
have width ${\cal O}(h^2)$ and are separated from each other by a
distance which is $\sim h$. We also recall from~\cite{Weinstein} and~\cite{Hitrik1} that 
\begeq
\label{4.3}
\Spec(P)\cap [\alpha,\beta] \subset \cup_{k=1}^{\infty} J_k(h).
\endeq

\noindent
When proving the bound (\ref{4.2}), we follow the suggestion of the remark
given after Theorem 3.2 in~\cite{Hitrik1}, which in turn is based on
the methods of Sections 2 and 3 of~\cite{Sjostrand1}. Thus, repeating the arguments
of the remark in~\cite{Hitrik1}, we find that there exists an elliptic selfadjoint $h$-pse\-udo\-diffe\-ren\-tial operator 
$A=\wrtext{Op}_h(e ^g)$, with $g=g(\Re z)\in S^{-1}(T^*M)$, such that   
$$
A ^{-1} (P+ih Q(z))A=P+ih\widehat{Q}+h^2 R(z),\quad R(z)\in \wrtext{Op}_h(S^0),
$$
where $\widehat{Q}=\wrtext{Op}_h(\hat{q})$ is selfadjoint, and $S^0\ni \hat{q}$ is such
that 
$$
\hat{q}=\hk{q}_T=\frac{1}{T}\int_0^T q\circ \exp(tH_p)\,dt \quad
\wrtext{on}\,\,\,p^{-1}(\Re z).
$$
Here $q=q(\Re z)=2a\sqrt{\Re z}$ and $T=T(\Re z)$ is the common period
of the closed $H_p$-trajectories in $p^{-1}(\Re z)$. As in the proof
of Proposition 3.1, we choose $g$ satisfying $H_pg=q-\hk{q}_T$ on
$p^{-1}(\Re z)$, and we may take for example,  
$$
g=\frac{1}{T} \int_0^T tq(\exp(tH_p))\,dt\quad\wrtext{on}\quad p^{-1}(\Re z).
$$
Since 
$$
A, A^{-1}={\cal O}(1): L^2 \rightarrow L^2,
$$
it suffices to show that the bound (\ref{4.2}) holds for the resolvent of 
$$
\widehat{{\cal P}}=P+ih\widehat{Q}+h^2 R(z).
$$

Now $H_p \hat{q}=0$ on $p^{-1}(\Re z)$, and we have 
$$
H_p\hat{q}=k(p-\Re z),
$$
with $k\in S^{-1}$. If we let $K\in \wrtext{Op}_h(S^{-1})$ have $k$ as
a principal symbol, then 
\begeq
\label{4.5}
ih[P,\widehat{Q}]=h^2K (P-\Re z)+h^3 S,
\endeq
where $S\in \wrtext{Op}_h(S^{0})$. Now, if $A$ and $B$ are bounded
selfadjoint operators, we have 
$$
\norm{(A+iB)u}^2=\norm{Au}^2+\norm{Bu}^2 +i([A,B]u,u)
$$
and using this identity together with (\ref{4.5}) we get, 
\begin{eqnarray}
\label{4.6}
& & 2\norm{(\widehat{{\cal P}}-z)u}^2 \geq
\norm{(P+ih\widehat{Q}-z)u}^2-{\cal O}(h^4)\norm{u}^2 \\ \nonumber 
& \geq & \norm{(P-\Re z)u}^2+
ih([P,\widehat{Q}]u,u) -{\cal O}(h^4)\norm{u}^2 \\ \nonumber
& \geq & \norm{(P-\Re z)u}^2-h\norm{[P,\widehat{Q}]u}\norm{u}-{\cal
O}(h^4)\norm{u}^2 \\ \nonumber
& \geq & \norm{(P-\Re z)u}^2-{\cal O}(h^2)\norm{(P-\Re z)u}\norm{u}
-{\cal O}(h^3)\norm{u}^2 \\ \nonumber
& \geq & \norm{(P-\Re z)u}^2 -{\cal O}(h)\norm{(P-\Re z)u}^2-{\cal
O}(h^3)\norm{u}^2 \\ \nonumber
& \geq & \frac{1}{2}\norm{(P-\Re z)u}^2-{\cal O}(h^3)\norm{u}^2, 
\end{eqnarray}
provided that $h$ is small enough. When $z$ is such that $\Im z={\cal
O}(h)$ and 
$$
\dist(\Re z,\cup_{k=1}^{\infty} J_k(h))\geq h/C,
$$
for some sufficiently large $C>0$, it follows from (\ref{4.3}) together with the spectral theorem that 
$$
\norm{(P-\Re z)u}\geq \frac{h}{C}\norm{u},
$$
and using this in (\ref{4.6}) we obtain, for $h$ small enough, 
$$
2\norm{(\widehat{{\cal P}}-z)u}^2\geq
\frac{h^2}{2C^2}\norm{u}^2-{\cal O}(h^3)\norm{u}^2 \geq
\frac{h^2}{4C^2}\norm{u}^2.
$$
This immediately implies (\ref{4.2}) with a new constant $C$. The proof is complete. 
\end{proof}

Repeating the arguments of Section 3 and using also (\ref{2.5.5}), we find that for $\tau$ as in
Proposition 4.1, as well as for all $\tau$ with $\Im \tau\geq
2\norm{a}_{\infty}+1$, we have 
\begeq
\label{4.7}
\norm{R(\tau)}_{{\cal L}(H^{\theta},L^2)}\leq \frac{{\cal
O}(1)}{(1+\abs{\tau})^{1+\theta/2}},\quad 0\leq \theta\leq 2. 
\endeq

We now come to prove the propagator expansion. When $k\in \nat$ is sufficiently large, we consider the positively
oriented contour $\gamma$ given by the following four line segments:
$\gamma_1=\{x-i; -k-\alpha/4-1/2\leq x \leq k+\alpha/4+1/2\}$, $\gamma_2=\{ k+\alpha/4+1/2+iy;
-1\leq y\leq 2\norm{a}_{\infty}+1\}$,
$\gamma_3=\{x+i(2\norm{a}_{\infty}+1); -k-\alpha/4-1/2\leq x\leq k+\alpha/4+1/2\}$, and
$\gamma_4=\{-k-\alpha/4-1/2+iy; -1\leq y\leq 2\norm{a}_{\infty}+1\}$. Let $f\in
H^{\theta}$, $\theta>0$. An
application of the residue theorem allows us to write 
\begin{eqnarray*}
& & \int_{\gamma_1} e^{it\tau} R(\tau)f\,d\tau=\sum_{\tau\in {\rm
Spec}({\cal A}),\, \abs{\Re \tau}< k+1/2+\alpha/4} e^{it\tau}p_{\tau}(t) f \\ \nonumber
& - & \int_{\gamma_2\cup \gamma_4} e^{it\tau}
R(\tau)f\,d\tau-\int_{\gamma_3} e^{it\tau} R(\tau)f\,d\tau.
\end{eqnarray*}
Here $p_{\tau}$ is a polynomial in $t$ with values in ${\cal
L}(L^2,H^s)$, for any $s\in \real$. The $L^2$-norm of the integral along the
segment $\gamma_2$ does not exceed 
$$
\int_{-1}^{2\norm{a}_{\infty}+1} e^{-ty} \frac{{\cal O}(1)\norm{f}_{H^{\theta}}}{(1+\abs{k})^{1+\theta/2}}\,dy,
$$
which tends to zero as $k\rightarrow \infty$. Similarly it follows
that the contribution coming from the segment $\gamma_4$ vanishes as
$k\rightarrow \infty$. When treating the contribution coming from the segment
$\gamma_3$, now with $k=\infty$, we notice that in view of
(\ref{4.7}) and since $\theta>0$, the contour $\gamma_3: \Im \tau=2\norm{a}_{\infty}+1$ can be replaced by $\Im \tau=C$ for any
$C\geq 2\norm{a}_{\infty}+1$. Letting $C\rightarrow \infty$, we conclude
that the corresponding contribution to the integral vanishes. We get
with $k_0\gg 1$, 
\begeq
\label{4.8}
U(t)f=\sum_{\tau \in {\rm Spec}({\cal A}),\,\, \abs{\Re
\tau}<k_0+\alpha/4} e^{it\tau} p_{\tau}(t)+\sum_{k=k_0}^{\infty} U_k(t)f,\quad \wrtext{when}\,\,t>0,\quad f\in H^{\theta},
\endeq
where
$$
U_k(t)=\sum_{\tau\in {\rm Spec}({\cal A}),\,\tau\in I_k}
e^{it\tau}p_{\tau,k}(t)+\sum_{\tau\in {\rm Spec}({\cal
A}),\,-\overline{\tau}\in I_k} e^{it\tau} p_{\tau,k}(t),
$$
and $p_{\tau,k}(t)$ is a polynomial in $t$ with values in ${\cal
L}(L^2,H^s)$, for any $s\in \real$. 
An application of (\ref{4.7}) together with the residue theorem shows that for each fixed $t$ we have
$$
\norm{U_k(t)}_{{\cal L}(H^{\theta},L^2)}={\cal O}(1)e^{t}
\frac{1}{\abs{k}^{1+\theta/2}}.
$$
The absolute convergence of the series (\ref{4.8}) in ${\cal
L}(H^{\theta},L^2)$, $\theta>0$, follows, and this concludes the proof of Theorem 1.4.

\section{Proof of Theorem 1.3}
\setcounter{equation}{0}
As in Section 1, we consider the operator 
$$
{\cal P}=P+ihQ(z), \quad P=-h^2 \Delta,\quad Q(z)=2a(x)\sqrt{z},\quad
z\in [\alpha,\beta]+i[-\gamma,\gamma],
$$
for $0<\alpha<1<\beta$, $\gamma>0$. In the following discussion, we shall assume that the $H_p$-flow possesses
a simple closed trajectory $\gamma:[0,T]\rightarrow p^{-1}(1)$, which is 
nondegenerate elliptic and such that 
\begeq
\label{5.1}
\pi(\gamma([0,T]))\cap \supp(a)=\emptyset.
\endeq
Here $\pi: T^*M\rightarrow M$ is the natural projection. Moreover, we
shall assume that the linearized Poincar\'e map $P_{\gamma}$ of
$\gamma$ is $N$-elementary, for some $N\geq 4$. We shall
establish the existence of a sequence of eigenvalues of the operator ${\cal A}$ 
converging rapidly to the real axis. 

Under these assumptions on $\gamma$, it is proved 
in~\cite{PopovCardoso} that there exist quasimodes for $P=P(h)$
associated with $\gamma$ and having polynomially small errors. To be
precise, it follows from the 
results of~\cite{PopovCardoso} that there exists a
set $H\subset (0,1]$ with $0\in \overline{H}$, the closure of $H$, 
and $0<a_0\leq a(h)\leq b(h)\leq b_0<\infty$, two functions on $H$, 
such that for any $h\in H$ there exists $m(h)\in \nat$, $E_j(h)\in
[a(h),b(h)]$, and $u_j(h)\in L^2$, $j=1,\ldots m(h)$, such that 
$$
\left(P(h)-E_j(h)\right)u_j(h)=R_0(h)\quad
\wrtext{in}\,\,\,L^2,\,\,j=1,\ldots\, ,m(h), 
$$ 
and $(u_j(h)|u_k(h))=\delta_{jk}+R_0(h)$, where $R_0(h)={\cal O}(h^{\infty})$. Furthermore, it follows from 
the construction of~\cite{PopovCardoso} that if $A\in {\rm Op}_h(S^0)$ is such that $WF_h(A)$ is 
contained in a small \neigh{} of $\gamma([0,T])$ and $A=I$ microlocally near $\gamma([0,T])$ then 
$$
(I-A)u_j(h)={\cal O}(h^{\infty}) \quad \wrtext{in}\,\,\,L^2,\,\,j=1,\ldots m(h).
$$
From the assumption (\ref{5.1}) we infer therefore that   
$$
({\cal P}-E_j(h))u_j(h)=R(h)\quad \wrtext{in}\,\,\,L^2,\quad j=1,\ldots m(h), 
$$
where $R(h)={\cal O}(h^{\infty})$. 

\Remark. The construction of~\cite{PopovCardoso} gives quasimodes of $P$
associated with a family of KAM tori of the Poincare
mapping, near the closed trajectory $\gamma$. For analytic manifolds, the construction
of~\cite{PopovCardoso} produces quasimodes with exponentially small
errors. We may also recall that under the stronger assumption of the
$N$-elementarity of $P_{\gamma}$ for every $N\in \nat$, the existence
of quasimodes with polynomially small errors can be established
also by the method of Gaussian beams---see~\cite{Ralston} and the references given there. 
While the construction of~\cite{PopovCardoso} gives quasimodes of
positive mass, given by the volume of the union of the KAM tori, it is known that the quasimodes obtained by the
method of Gaussian beams have mass zero---see~\cite{CdV}.

\vskip 4mm
For future reference, we notice that for $z\in
[\alpha,\beta]+i[-\gamma,\gamma]$, we have $\Im ({\cal P}u|u)\geq 0$, and therefore
\begeq
\label{5.2}
\norm{({\cal P}-z)^{-1}}_{{\cal L}(L^2,L^2)}\leq \frac{1}{\abs{\Im z}}, \quad \Im z<0.
\endeq

The following result implies Theorem 1.3. 

\begin{theo}
For any positive function $S(h)$ with $De^{-D/h}<S(h)={\cal O}(h^{\infty})$,
for some $D>0$, and such that $S(h)\gg R(h)$, and every $k\in \nat$, there exists $h(S,k)>0$ such that for
$H\ni h\leq h(S,k)$, the spectrum of ${\cal P}={\cal P}(h)$ intersects the set 
\begeq
\label{5.3}
[E_j(h)-6h^k,E_j(h)+6h^k]+i[0,2S(h)h^{-n-1}],\quad j=1,\ldots m(h).
\endeq
\end{theo}

When proving Theorem 5.1, we shall start by establishing the following
result, which gives an a priori exponential estimate of the resolvent of ${\cal P}$, in a
set of size $\sim h$, outside the union of small neigbourhoods of the eigenvalues. 

\begin{prop}
Let $E\in [\alpha+1/C_0,\beta-1/C_0]$ for some $C_0>1$. For every $L>0$ there exists $A>0$ and $h_0>0$ such that for $0<h<h_0$
we have 
\begeq
\label{5.4}
\norm{({\cal P}-z)^{-1}}_{{\cal L}(L^2,L^2)}\leq \frac{A}{h} \exp \left(A h^{1-n}
  \log\frac{h}{g(h)}\right),
\endeq
for $\abs{z-E}<Lh$, $\dist(z, {\rm Spec}({\cal P}))\geq g(h)$. Here
$0<g(h)\ll h$.  
\end{prop}
\begin{proof}
We shall use the idea of~\cite{Markus} and~\cite{Sjostrand1} of exploiting a finite rank perturbation to 
create a gap in the spectrum of $P$ around $E$. Recall from~\cite{Sjostrand1} that for any $C>0$ there exists a 
selfadjoint operator $\widetilde{P}=P+K$, with the same domain as $P$, such that 
$$
[E-Ch,E+Ch]\cap \Spec(\widetilde{P})=0,
$$
and
$$
\norm{K}_{{\cal L}(L^2,L^2)}\leq Ch,\quad \norm{K}_{\mathrm{tr}}\leq \widetilde{C}(C)h^{2-n}. 
$$
Here $\norm{\cdot}_{\mathrm{tr}}$ denotes the trace class norm. We then introduce 
$$
\widetilde{{\cal P}}=\widetilde{P}+ihQ(z), 
$$
and writing 
$$
\widetilde{{\cal P}}-z=(\widetilde{P}-z)(I+ih(\widetilde{P}-z)^{-1}Q),
$$
we see that 
\begeq
\label{5.4.5}
(\widetilde{{\cal P}}-z)^{-1}=\frac{{\cal O}(1)}{h}
\quad \wrtext{in}\quad {\cal L}(L^2,L^2), 
\endeq
when $z$ is such that $\abs{\Re z-E}\leq Lh$, provided that $C$ is
large enough. Now we write
\begeq
\label{5.5}
{\cal P}-z=({\cal \widetilde{P}}-z)(I-({\cal \widetilde{P}}-z)^{-1}K),
\endeq
and let us introduce 
$$
D(z)=\det(I-({\cal \widetilde{P}}-z)^{-1}K),
$$
which is a holomorphic function in the open disc $D(E,Lh)=\{z\in
\comp; \abs{z-E}<Lh\}$, whose zeros are precisely the eigenvalues of
${\cal P}$ in this set. Since the trace class norm of 
$({\cal \widetilde{P}}-z)^{-1}K$ is ${\cal O}(1)h^{1-n}$, we get 
\begeq
\label{5.5.5}
\abs{D(z)}\leq \exp({\cal O}(1)h^{1-n}),\quad z\in D(E,Lh).
\endeq
An application of Theorem 5.1 from Chapter 5 of~\cite{GoKr} shows that
when $\abs{\Re z-E}\leq Lh$, 
\begeq
\label{5.6}
\norm{(I-({\cal \widetilde{P}}-z)^{-1}K)^{-1}}_{{\cal L}(L^2,L^2)} 
\leq \frac{\det(I+\abs{({\cal
\widetilde{P}}-z)^{-1}K})}{\abs{D(z)}}\leq \frac{e^{{\cal
O}(1)h^{1-n}}}{\abs{D(z)}}, 
\endeq
and when deriving bounds on the resolvent of ${\cal P}$, we shall
therefore have to
estimate $\abs{D(z)}$ from below, away from its zeros. In doing so we
remark that when $\abs{z-E}<Lh$ and $\abs{\Im z}\geq {\cal O}(1)h$ then 
$$
({\cal P}-z)^{-1}=\frac{{\cal O}(1)}{h}: L^2 \rightarrow L^2,
$$
and from (\ref{5.5}) we get that in this set, 
$$
\left(I-({\cal \widetilde P}-z)^{-1}K\right)^{-1}=({\cal P}-z)^{-1} ({\cal \widetilde P}-z)=I+({\cal P}-z)^{-1}K,
$$
with 
$$
\norm{({\cal P}-z)^{-1}K}_{\mathrm{tr}}\leq {\cal O}(1) h^{1-n}.
$$
It follows that 
\begeq
\label{5.6.5}
\abs{D(z)}\geq e^{-{\cal O}(1)h^{1-n}}, \quad z\in D(E,Lh),\,\,\,\abs{\Im z}\geq {\cal O}(1)h.
\endeq
Let now $z_1,\ldots z_N$ be the zeros of $D(z)$ in $D(E,Lh)$. From~\cite{Sjostrand1} we recall that $N={\cal O}(1) h^{1-n}$. We 
factorize
$$
D(z)=G(z)B(z), \quad z\in D(E,Lh), 
$$
where $G$ is holomorphic and non-vanishing in $D(E,Lh)$, and 
$$
B(z)=\prod_{j=1}^N b_{z_j}(z)
$$
is the Blaschke product on $D(E,Lh)$, so that
$$
b_{z_j}(z)=\hat{b}_{\frac{z_j-E}{Lh}}\left(\frac{z-E}{Lh}\right)=
\frac{z-z_j}{Lh\left(1-\left(\frac{z-E}{Lh}\right)\left(\frac{\overline{z_j}-E}{Lh}\right)\right)}.
$$
Here
$$
\hat{b}_w(z)=\frac{z-w}{1-\overline{w}z}, \quad \abs{z}\leq 1, \,\,
\abs{w}<1, 
$$
stands for the standard Blaschke factor for the unit disc. An
application of basic factorization theorems in the Banach algebra of
bounded holomorphic functions in a disc, or simply using the maximum
principle together with (\ref{5.5.5}), allows us to conclude that 
\begeq
\label{5.7}
\abs{G(z)}\leq e^{{\cal O}(1)h^{1-n}} \quad \wrtext{in}\,\,\,D(E,Lh).
\endeq
Taking (\ref{5.6.5}) into account we find that for $\abs{\Im z}\geq {\cal O}(1)h$,  
$$
\abs{G(z)}\geq e^{-{\cal O}(1) h^{1-n}},
$$
An application of the Harnack inequality to the harmonic function
$\log\abs{G(z)}$, in the form given in Chapter 1 of~\cite{Markus}, shows that 
$$
\abs{\log \abs{G(z)}}\leq {\cal O}(1) h^{1-n}, \quad \wrtext{in}\,\,D(E,Lh),
$$
after decreasing $L$ slightly. It is easy to see that 
$$
\abs{B(z)}\geq e^{-{\cal O}(1)h^{1-n} \log{\frac{h}{g(h)}}},\quad z\in
D(E,Lh),\,\,\,\dist(z,{\rm Spec}({\cal P}))\geq g(h),
$$
and therefore for such $z$'s we get, 
\begeq
\label{5.8}
\abs{D(z)}\geq e^{-{\cal O}(1)h^{1-n} \log{\frac{h}{g(h)}}}. 
\endeq
Combining (\ref{5.5}) and (\ref{5.4.5}) with (\ref{5.6}) and
(\ref{5.8}), we complete the proof of (\ref{5.4}). 
\end{proof}

\Remark. The resolvent bound of Proposition 5.2 together with its
proof is closely related to the abstract resolvent estimates for weak
non-selfadjoint perturbations of selfadjoint operators, developed in Chapter 1 of~\cite{Markus}. See
also~\cite{BurqZworski} and~\cite{DSZ} for further references and
applications of similar bounds in the theory of resonances and
semiclassical pseudospectral theory.

\Remark. For the following argument, it would have been sufficient to
derive an exponential bound on $({\cal P}-z)^{-1}$ in a bounded
$h$-independent set, rather than in a set of size $\sim h$. Indeed, if
$\Omega\subset \subset \{z\in \comp; \Re z>0\}$ is a bounded open $h$-independent
neighbourhood of some $E>0$, then a straightforward adaptation of the proof
of Proposition 5.2 shows that for all $h$ small enough and $z\in
\Omega$, we have 
\begeq
\label{5.8.5}
\norm{({\cal P}-z)^{-1}}_{{\cal L}(L^2,L^2)}\leq {\cal O}(1)
\exp({\cal O}(1) h^{-n} \log{(1/g(h))}),\quad \dist(z,{\rm Spec}({\cal P}))\geq g(h)>0.
\endeq
See also Lemma 6.1 of~\cite{DSZ}. The exponent $h^{-n}$ in
(\ref{5.8.5}) reflects the fact that the number of eigenvalues of ${\cal P}$ in $\Omega$ is ${\cal
O}(h^{-n})$---see~\cite{Sjostrand1} for the much more precise
result giving the Weyl asymptotics. It seemed to us nevertheless worthwhile,
and indeed, of independent interest, to obtain a sharper result in the
smaller set, whose size is dictated by the strength of the
non-selfadjoint perturbation. Notice also that in the one-dimensional case,
$n=1$, Proposition 5.2 gives a polynomial bound on the resolvent of 
${\cal P}$. It follows that in this case, only eigenvalues can produce a
super-polynomial growth of the resolvent.   

\vskip 2mm
\noindent
We now come to prove Theorem 5.1. When doing so, we follow the
argument of~\cite{TZ1} and argue by contradiction. Theorem 5.1 is then
obtained by applying the argument of~\cite{TZ1} as it stands, since an application of the semiclassical maximum
principle, Lemma 2 of~\cite{TZ1}, to $({\cal P}-z)^{-1}$ is now
legitimate, in view of 
Proposition 5.2, (\ref{5.8.5}), and the bound (\ref{5.2}). The proof of
Theorem 5.1 is complete. 

\Remark. Following~\cite{Stefanov}, it is possible to sharpen the
statement of Theorem 5.1 and to estimate the number of eigenvalues of
${\cal P}$ in the set
\begeq
\label{5.10}
[a(h)-{\cal O}(1)h^k, b(h)+{\cal O}(1)h^k]+i[0,2S(h)h^{-n-1}],\quad
k\in \nat. 
\endeq
Below we shall indicate how to do so, by adapting the argument of~\cite{Stefanov}
to the present setting. To this end, let us introduce the operator
\begeq
\label{5.11}
{\cal A}(h)= \left( \begin{array}{cc}
   0 & 1 \\
  -h^2 \Delta & 2iah
\end{array} \right)=\left( \begin{array}{cc}
   1 & 0 \\
  0 & h
\end{array} \right) h{\cal A}\left( \begin{array}{cc}
   1 & 0 \\
  0 & h
\end{array} \right)^{-1}: H^1\times L^2 \rightarrow H^1\times L^2.
\endeq
It is then true that $z\in {\rm Spec}({\cal P})$ if and only if
$\lambda=\sqrt{z}\in {\rm Spec}({\cal A}(h))$. Using the factorization
in (\ref{5.11}) together with Proposition 2.1, we verify that 
$$
\norm{({\cal A}(h)-\lambda)^{-1}}_{{\cal L}({\cal H},{\cal H})}\leq
\frac{{\cal O}(1)}{h\abs{\Im \lambda}},\quad \Im \lambda<0.
$$
Here ${\cal H}=H^1\times L^2$, and we equip $H^1$ with the
corresponding $h$-dependent norm. It is also easy to see, using
Proposition 5.2, that the resolvent of ${\cal A}(h)$
enjoys an exponential type bound of the same form as $({\cal
P}-z)^{-1}$, when
away from the eigenvalues. If we now use the quasimodes for $P$ and set
\begeq
\label{5.12}
U_j(h)=\left( \begin{array}{cc}
u_j(h) \\
E^{1/2}_j(h)u_j(h)
\end{array}\right)\in {\cal H},\quad j=1,\ldots m(h),
\endeq
then a simple computation shows that 
$$
({\cal A}(h)-E^{1/2}_j(h))U_j(h)=R(h)\quad \wrtext{in}\,\,\,{\cal H},
$$
and $(U_j(h)|U_k(h))_{{\cal H}}=(1+E_j(h)+E_k(h))\delta_{jk}+R(h)$,
where $R(h)={\cal O}(h^{\infty})$. Here we have used that the
modes $u_j(h)$ are orthonormal in $L^2$ and ortogonal in $H^1$, modulo
an error which is ${\cal O}(h^{\infty})$. 

An inspection of~\cite{Stefanov} combined with the
preceding remarks shows that the coun\-ting argument of Theorem 1
of~\cite{Stefanov} applies to the operator-valued function $h({\cal
A}(h)-\lambda)^{-1}$ without any change, and demonstrates that for any positive function $S(h)$
such that $\max(De^{-D/h},h^{-n-2}R(h))\leq S(h)={\cal
O}(h^{\infty})$, $D>0$, and any $k\in \nat$, the rectangle 
$$
[a(h)^{1/2}-6h^k,b(h)+6h^k]+i[0,2S(h)h^{-n-1}]
$$
contains $\geq m(h)$ eigenvalues of ${\cal A}(h)$, for $h$ small
enough. It follows that the
set (\ref{5.10}) contains at least $m(h)$ eigenvalues of ${\cal P}$.


\end{document}